\documentclass[12pt,draft]{amsart}
\usepackage{amsmath,amsthm,latexsym,amscd,amsbsy,amssymb,amsfonts,
euscript}

\newcommand{\GB}[1]{G^{\B }_{#1}}
\newcommand{\GBP}[1]{G^{\B '}_{#1}}
\renewcommand{\AE}{\mathop{\rm AE}\nolimits}
\newcommand{\ANE}{\mathop{\rm ANE}\nolimits}
\newcommand{\dist}{\mathop{\rm dist}\nolimits}
\newcommand{\St}{\mathop{\rm St}\nolimits}
\newcommand{\dimm}{\mathop{\rm dim}\nolimits}

\newcommand{\A}{{\mathfrak A}}
\newcommand{\B}{{\mathfrak B}}
\newcommand{\G}{{\EuScript G}}
\newcommand{\F}{{\EuScript F}}
\newcommand{\U}{{\EuScript U}}
\newcommand{\V}{{\EuScript V}}
\newcommand{\ed}{\mathop{\rm ed}\nolimits}
\newcommand{\diam}{\mathop{\rm diam}\nolimits}
\newcommand{\cit}[1]{{\rm {(\cite{#1})}.}}
\newcommand{\cov}{\mathop{\rm cov}\nolimits}
\newcommand{\ord}{\mathop{\rm ord}\nolimits}
\newcommand{\e}{\varepsilon}
\newcommand{\al}{\alpha}
\renewcommand{\d}{\delta}

\newcommand{\be}{\beta}

\newcommand{\x}{\xi}

\newcommand{\m}{\mu}
\newcommand{\n}{\nu}
\renewcommand{\l}{\lambda}
\newcommand{\emp}{\varnothing}

\newtheorem*{MainLemma}{\indent Main Lemma}
\newtheorem*{Tm}{Theorem}
\newtheorem{thm}{Theorem}[section]
\newtheorem{cor}[thm]{Corollary}
\newtheorem{lem}[thm]{Lemma}
\newtheorem{pro}[thm]{Proposition}

\theoremstyle{definition}
\newtheorem{dfn}{Definition}[section]
\theoremstyle{remark}
\newtheorem{rem}{Remark}[section]

\chardef\bslash=`\\ 

\makeatletter
\def\verbatim{\interlinepenalty\@M \@verbatim
  \leftskip\@totalleftmargin\advance\leftskip2pc
  \frenchspacing\@vobeyspaces \@xverbatim}
\makeatother
\numberwithin{equation}{section}

\begin{document}


\title[Continuous selections with respect to extension dimension]
{Continuous selections with respect to extension dimension}
\author{A.~V.~Karasev}
\address{Department of Mathematics and Statistics,
University of Saskatche\-wan,
McLean Hall, 106 Wiggins Road, Saskatoon, SK, S7N 5E6,
Canada}
\email{karasev@math.usask.ca}

\keywords{Continuous selections, extension dimension,
dimension}
\subjclass{Primary: 54C65; Secondary: 54C20, 54F45}


\begin{abstract}
Let $L$ be finite CW complex. By $[L]$ we denote extension type of $L$.
The following generalization of Michael's selection theorem is proved:

\begin{Tm}
Consider $[L]\le [S^n]$. Let $F$ be a lower semi-continuous map between
polish space $X$ and metrizable compactum $Y$, such that $\ed (X) \le
[L]$, $F$ is equi-$LC^{[L]}$ collection and $F(x)\in \AE ([L])$ for any
$x \in X$. Let $A$ be a closed subset of $X$ such that there exists
a continuous selection $f:A \to Y$ of $F|_A$.
Then $F$ admits a continuous selection $\bar
f$ which extends $f$.
\end{Tm}

\end{abstract}

\maketitle
\markboth{A.~V.~Karasev}{Continuous selections with respect to
extension dimension}


\section{Introduction}

The following Michael's selection theorem is well-known (see
\cite{mchl} for details):

\begin{thm} \label{mhl}
Let $X$ be a paracompact space, $A\subseteq X$ a closed subspace of
$X$ with  $\dimm _X (X - A) \le n+1$,
$Y$ a complete metric space, $F$ an equi-$LC^n \bigcap C^n$
collection and $F:X \to Y$ lower semi-continuous map.
Then every selection for $F|_A$
can be extended to a selection for $F$.
\end{thm}

Theorem~\ref{mhl} deals with usual Lebesgue dimension and concerns the notion
of absoulute extensor in dimension $n$.
The purpose of the present paper is to obtain a natural generalization
of Michael's theorem in the case of extension dimension.

\section{Preliminaries}
     In this part we introduce notions of {\it extension types of complexes,
extension dimension, absolute extensors modulo a complex, $[L]$-homotopy}  and
{\it equi-$LC^{[L]}$ collections}.
All spaces are polish, all complexes are countable finitely-dominated
$CW$ complexes.
For more details related to extension dimension
see \cite{ch}.

For spaces $X$ and $L$, the notation $L \in \AE (X)$ means, that every map
$f:A \to L$, defined on a closed subspace $A$ of $X$, admits an
extension $\bar f$ over $X$.

Let $L$ and $K$ be complexes. We say (see \cite{ch}) that $L \le K$ if
for each space $X$ from $L \in \AE (X)$ follows $K \in \AE (X)$. Equivalence
classes of complexes with respect to this relation are called {\it
extension types}. By $[L]$ we denote extension type of $L$.

\begin{dfn} \cit {ch} The extension dimension of a space $X$ is
extension type $\ed (X)$ such that
$\ed (X) = \min \{ [L] : L \in \AE  (X) \}$.
\end{dfn}

Observe, that if $[L] \le [S^n]$ and $\ed (X) \le [L]$, then $\dimm X
\le n$.

\begin{dfn}\cit {ch}
We say that a space $X$ is {\it an absolute extensor modulo}
$L$ (shortly $X$ is  $\AE ([L])$)
and write $X\in \AE ([L])$ if $X\in \AE (Y)$ for each space
$Y$  with $\ed (X) \le [L]$.
\end{dfn}

We will widely use the following proposition:

\begin{pro} \label{sum} \cit {ch}
Let $X$ be a polish space such that $\ed (X) \le [L]$ and $Y \in
\AE ([L])$. Then $L \in \AE (X')$ for any $X' \subseteq X$.
\end{pro}

Follow \cite {ch} give definition of $[L]$-{\it homotopy} and
$u$-$[L]$-homotopy:
\begin{dfn}
Two maps $f_0$, $f_1: X \to Y$ are said to be $[L]$-homotopic (notation:
$f_0 \stackrel{[L]}{\simeq} f_1$) if for any map
$h: Z \to X \times [0,1]$, where $Z$ is a space with $\ed (Z) \le [L]$,
the composition $(f_0 \oplus f_1) h|_{h^{-1}(X \times \{ 0, 1 \} )} :
h^{-1} (X \times \{ 0,1 \}) \to Y$ admits an extension $H: Z \to Y$.
If, in addition, we are given  $\U \in \cov (Y)$
and $H$ can be choosen so that the
collection $\{ H(h^{-1} (\{ x \} \times [0,1])): x \in X \}$
refines $\U $, we say, that $f_0$ and $f_1$ are $\U$-$[L]$-homotopic, and
write $f_0 \stackrel{\U -[L]}{\simeq} f_1$.
\end{dfn}

It is clear, that if $f_0$, $f_1$ are $\U $-$[L]$-homotopic for some
$[L]$ then these maps are $\U $-close.

Let us observe (see \cite {ch}) that $\AE ([L])$-spaces have the following
important property:

\begin{pro}\label{homot}
Let $Y$ be a Polish $\AE ([L])$-space. Then for each $\U \in \cov (Y)$
there exists $\V \in \cov (Y)$ refining $\U$, such that for any space
$X$ with $\ed (X) \le [L]$, any closed subspace $A \subseteq X$ and any
two $\V$-close maps $f, g: A \to Y$ from  existance of extension $\bar
f$ of $f$ over $X$ follows existance $\bar g$  which is $\U $-close to
$\bar f$ and extends $g$ over
$X$.
\end{pro}

\begin{cor}\label{c1}
Let $Y$ be a copmact $\AE ([L])$-space.\\
Then for each $\e > 0$
there exists $\d = \d (\e ) > 0$  such that for any space
$X$ with $\ed (X) \le [L]$, any closed subspace $A \subseteq X$ and any
two $\d$-close maps $f, g: A \to Y$ from  existance of extension $\bar
f$ of $f$ over $X$ follows existance $\bar g$  which is $\e $-close to
$\bar f$ and extends $g$ over
$X$.
\end{cor}

 From just mentioned fact one can easely obtain the following:

\begin{pro}\label{hext}
Let Y be a metrizable $\AE ([L])$ compactum. Then for each $\e > 0$
there exists
$\d = \d (\e ) > 0$  such that for any space
$X$ with $\ed (X) \le [L]$,
any closed subspace $A \subseteq X$ and any
 map $f: A \to Y$  such that
$\diam (f(A)) \le \d $ there exists $\bar f: X \to Y$  extending $f$
 such that $ \diam (\bar f (X) \le \e$.
\end{pro}

The last proposition allows us to introduce in the natural way the notion of {\it
equi-}$LC^{[L]}$ {\it collection}.
\begin{dfn}\label{eq}
Collection $\F = \{ F_{\al} : \al \in \A \}$ of closed sabsets of
compact space $Y$ is said to be
equi-$LC^{[L]}$ if for any $\e > 0$ there exists $\d >0$ such that for
each $\al \in \A$, each Polish space $Z$ with $\ed (Z) \le [L]$, each
closed subset $A \subseteq Z$ and a map $f:A \to F_{\al}$ such that
$\diam (f(A)) < \d$ there exists an extension $\bar f: Z \to F_{\al}$
of $f$ over $Z$ such that $\diam f(Z) < \e$.
\end{dfn}

\section{Selection theorem}
Let us recall that a many-valued map  $F:X \to Y$ is said to be lower
semi-continuous (shortly l.s.c.) if $F(x)$ is closed subset of $Y$
for any\\ $x \in X$ and for any open $U \subseteq Y$ the set
$F^{-1} (U) = \{ x \in X: F(x) \bigcap U \ne \emp \}$ is open in $X$.

We are ready now to formulate our main result.

\begin{thm}\label{mt1}
Let $Y$ be a metrizable  compactum, $X$ a Polish space with
$\ed (X) \le [L]$ and $F:X \to Y$ a
l.s.c. map such that collection $\F = \{ F(x): x \in X \}$ is
equi-$LC^{[L]}$ and $F(x) \in \AE ([L])$ for each $x\in X$. Let $A \subseteq X$
be a closed subset. Then  any
selection $f:A \to Y$ of $F|_A$  can be extended to selection $\bar f:X \to Y$.
\end{thm}

\begin{rem}
Important difference between the further proof of Theorem \ref{mt1}
and consideration of \cite{mchl} consists in the fact that we cannot
apply technique of \cite{mchl} involving maps into nerves of covering.
We have to directly extend maps over open subspaces of $X$ and hence
we need to use Proposition \ref{sum}.
Therefore proof presented in this text cannot be directly generalized
on the case when $X$ is paracompact space.

Further, we have no characterization of absolute extensors modulo $[L]$
in terms of maps of spheres.
It makes Proposition \ref{hext} and in turn compactness of $Y$
essential for our consideration.

It is also necessary to point out, that in \cite{ch} $[L]$-dimensional
analogies of $n$-dimensional spheres were introduced, namely, a compact
spaces $S^n_{[L]}$, which are $\ANE ([L])$ and admit $[L]$-invertable
and approximately $[L]$-soft mappings onto $n$-dimensional sphere (see
\cite{ch} for necessary definition).
Additionaly, these spaces are proved to be $[L]$-universal for compact
spaces.
This fact, it would seem, allows to introduce the notion of
equi-$LC^{[L]}$ families using characterization in terms of mappings of
$S^n_{[L]}$ which were closer to original definition in \cite{mchl}, and
generalize the theorem on the case of non-compact $Y$.

Unfortunately, as it already has been mentioned above,
our proof involve extansions of maps over open subspaces of $X$ which
are non-compact. Therefore we cannot use universality of $S^n_{[L]}$.

\end{rem}

Simillar \cite{mchl}, we accomplish the proof of this theorem
consequently reducing it to other assertion.
Using arguments of  \cite{mchl}, one can easely observe, that Theorem
\ref{mt1} is equivalent to the following

\begin{thm}\label{mt2}
Let $Y = Q$ be Hilbert cube, $X$ a Polish space with
$\ed (X) \le [L]$ and $F:X \to Q$ an
l.s.c. map such that collection $\F = \{F(x): x \in X \}$ is
equi-$LC^{[L]}$ and $F(x) \in \AE ([L])$ for each $x\in X$.
Then $F$ admits selection  $f:X \to Y$.
\end{thm}

Let $B \subseteq Y$. By $O_{\e} (B)$ we denote $\e$-nighbourhood of $B$
in Y.

Finally, let us reduce Theorem~\ref{mt2} to the following lemma:

\begin{MainLemma}
Let $X$, $Y$ and $F$ be the same as in Theorem~\ref{mt2}. Then

{
\renewcommand{\theenumi}{{\bf \alph{enumi}}}
\begin{enumerate}
\item  For any $\m > 0$ there exists $g:X \to Y$, which is $\m$-close
to $F$.
\item  For any $\e > 0$ there exists $\d = \d (\e ) > 0$
with the following\\ property:
for each $f: X \to Y$ such that
$f$ is $\d$-close to $F$  and for each $\m > 0$
there exists $g: X \to Y$
such that $g$ is $\e$-close to $f$ and $\m$-close to $F$.
\end{enumerate}
}
\end{MainLemma}

Let us prove that

\begin{pro}
Main Lemma implies Theorem~\ref{mt2}.
\end{pro}

\begin{proof}
Consider a sequence $ \e _n = \frac{1}{2^n}$, $n \ge 0$. Using Main
Lemma construct corresponding sequences of $\{ \d _n < \e _n \}$, where
$\d _n = \d _n (\e _n)$ and
$\{ f_n \}$
such that $f_n$ is $\e _n$-close to $f_{n-1}$, $n \ge 1$ and $\d _n$ close
to $F$ for every $n$.    Then $f_n$ is uniformly Cauchy. Since $Y$ is
metrizable compactum (actually we assume that $Y$ is Hilbert cube),
there exists continuous $f = \lim\limits_{n \to \infty} f_n $. Obviously, $f$ is selection
of $F$.
\end{proof}

\section{Covers of special type}\label{cov}

Let us introduce notations and definitions
which are necessary to prove Main Lemma.

Since this point and up to the end of the text we assume that $L$ is a
complex such that $[L] \le [S^n]$, $X$ is
a Polish space with $\ed (X) \le [L]$ (and therefore with $\dimm Xf
\le n$), $Y$ is Hilbert cube (actually we need only the property $Y \in
\AE $ and compactness of $Y$)
and $F:X \to Y$ as in formulation of Main Lemma.

\begin{dfn}\label{can} Let $\U \in \cov (X)$. Then $\V \in \cov (X)$
is said to be {\it a canonical refinment}
for $\U$ if $\V$ satisfies the following conditions:

\begin{enumerate}
\item \label{st1} $\V$ is star-refinment of $\U$.
\item \label{stf1} $\V$ is star-finite.
\item \label{ord1} Order of $\V$ is $\le n+1$.
\item \label{ir2} $\V$ is irreducible, i.e. for any $V \in \V$ collection $\V \backslash \{
V \} $ is
not a cover of $X$.
\end{enumerate}
\end{dfn}

Observe, that canonical refinment exists for any $\U \in \cov (X)$ (see
\cite{eng} for details).

For any $\e > 0$ let $U_x = \{x' \in X: F(x) \subseteq  O_{\e}
(F(x'))\}$.   Since $F$ is l.s.c., $U_x$ is open for any $x\in X$
\cite{mchl}.

\begin{dfn}\label{canF} Let $\U \in \cov (X)$.
Then we say that $\V \in \cov (X)$
is  {\it a canonical refinment
for} $\U$ {\it with respect to} $F$ {\it and} $\e $
(notation: $\V = \V (U, F, \e)$, if it satisfies the following
conditionsf:

\begin{enumerate}
\item \label{st2} $\V$ is star-refinment of $\U$.
\item \label{ex2} For any $V \in \V$ there exists $x(V) \in V$ such that $ V
\subseteq U_{x(V)}$.
\end{enumerate}
\end{dfn}

It is easy to see, that for any $\U \in \cov (X)$  and $\e >0$ there exists a
canonical refinment with respect to $F$ and $\e$.

Let $\U = \{ U_{\al}: \al \in \A \}$ be a star-finite irreducible cover of $X$
having oder $\le n+1$.

Let $F_k = \{ x \in X : \ord _{\U} x \le k \}$, $k = 1 \ldots n+1 $.
Observe, that

{
\renewcommand{\labelenumi}{F\theenumi .}
\begin{enumerate}
\item \label{f1} $F_k$ is closed for any $k$.
\item \label{f2} $F_k \subseteq F_{k+1}$ for any $k = 1 \ldots n$.
\item \label{f3} $X = F_{n+1}$.
\end{enumerate}
}

Further, for each $k = 1 \ldots n+1$ and $\B \subseteq \A$ such that $|\B | = k$
let $G^{\B}_k = F_k \bigcap (\bigcap \{ U_{\al} : \al \in \B \} )$.

Notice that generally speaking, $G^{\B}_k$ may be empty or non-closed.

Obviously, family $\G _k= \{ G^{\B}_k : |\B | =k \} $ has the following
properties:

{
\renewcommand{\labelenumi}{G\theenumi .}
\begin{enumerate}
\item \label{G1} $\{ G^{\{ \al \}}_1 : \al \in \A \}$ are closed,
pairwise disjoint and non-empty subsets of $X$.
\item \label{G2} $\G _k$ is discrete in itself.
\item \label{G3} $F_{k+1} = (\bigcup \G _{k+1}) \bigcup F_{k}$
\item \label{G4} For each non-empty $\GB{k+1}$,
$\bigcup\{ \GBP{|\B '|}: \B ' \subsetneq \B \}
\supseteq \overline{\GB{k+1}}\bigcap F_k \ne \emp$
\end{enumerate}
}

We will use these consideration as well as introduced notations
in all the remaining text.

\section{Technical lemmas}

The following two lemmas we need to complete the proof are analogies of lemmas containing in Appendix
of \cite{mchl}.

\begin{lem}\label{l1}
Let $B$ be closed subset of $Y$, such that $B \in \AE ([L])$. Then for
any $\e > 0$ there exists $\d = \d (\e ) > 0$ such that for every Polish
space $X$ with $\ed (X) \le [L]$ and for each $f:X \to O_{\d}(B)$ there
exists $g: X \to B$ such that $g$ is $\e$-close to $f$.
\end{lem}

\begin{proof}
Construct a sequence $\{ \d _k \}^{n+1}_{k=1}$ such that:

{
\renewcommand{\labelenumi}{\theenumi ${\rm \d}$.}
\begin{enumerate}
\item \label{d1} $\d _{n+1} = \e$.
\item \label{d2} $\d _k < \d _{k+1}$.
\item \label{d3} Pair $(\frac{\d_{k+1}}{3}, \d_k) $ satisfies condition of
Proposition~\ref{hext}.
\end{enumerate}
}

Let $\d = \frac{\d _1}{6}$.
Check that pair $(\e , \d)$ satisfies requirments of lemma. Consider $f:
X \to O_{\d } (B)$.

Let $\EuScript O = \{ O_y =  O_{\d } (y) : y \in B\}$ (since $B$ is compact,
we may choose finite refinment of $\EuScript O$, but it is not essential
for further consideration).
Let $V_y = f^{-1} (O_y)$ and $\V = \{ V_y: y \in B \}$.
Consider $\U \in \cov (X)$, which is canonical refinement of
$\V$ in the sense of Definition~\ref{can}. Let $\U = \{ U_{\al} : \al
\in \A \}$.
Using property~\ref{st1} of canonical refinement, for each $\al \in \A$
find $y_{\al}$ such that $\St _\U U_{\al} \subseteq V_{y_{\al}}$. Notice
that generally speaking, $y_{\al}$ and $y_{\be}$ may coincide for $\al
\ne \be$.
Finally, consider related to $\U$ sets $F_k$ and
$G^{\B}_k$, introduced in Section~\ref{cov}.

We are ready now to construct map $g$.

Using induction by $k = 1 \ldots n+1$ construct a sequence
of map $\{ g_k \} ^{n+1}_{k=1}$ such that:

{
\renewcommand{\labelenumi}{\theenumi {\rm g}.}
\begin{enumerate}
\item \label{1g} $g_k: F_k \to B$.
\item \label{2g} $g_{k+1}|_{F_k} = g_k$.
\item \label{3g} $f|_{F_k}$ is $\e$-close to $g_k$.
\item \label{4g} $g_k(F_k \bigcap U_{\al}) \subseteq
O_{\d _k/3}(y_{\al})$ for each $\al \in \A$.
\end{enumerate}
}

For $k = 1$ define $g_1$ letting $g_1|_{G^{\{\al \}}_1} \equiv y(\al)$.
Observe, that $g_1$ is defined correctly and continuous on $F_1$
(see properties G\ref{G1}--G\ref{G4} on the page \pageref{G1}).
By our choice of $\{ y_{\al}\}$, $g_1$ satisfies conditions
\ref{1g}g--\ref{4g}g.

Assuming that $g_k$ has been already constructed, let us construct
$g_{k+1}$.

To accomplish this it is enough (see G\ref{G3} on the page  \pageref{G3})
to define $g_{k+1}$ on each non-emty $G^{\B }_{k+1}$ for each $\B \subseteq \A$
such that $|\B | = k+1$.

Fix $\B$ such that $|\B | = k+1$. Consider $Z = \bigcup\limits_{\al \in
\B} U_{\al}$.

Let $Z' = Z \bigcap F_k$. Obviously, $Z'$ is closed subset
of $Z$. Since $U_{\al} \bigcap U_{\be} \ne \emp$ for $\al , \be \in \B$
(recall that we consider $G^{\B }_{k+1} \ne \emp$) by our choice of
$y_{\al}$ we have $\dist (y_{\al},y_{\be}) < \d$. Therefore, by
property~\ref{4g}g we conclude that $\diam Z' < 2\d +2(\d _k /3) < \d
_k$. Hence by our choice of $\{\d_ k\}$, map $g_k$ has an extension
$\bar g_k : Z \to B$ such that $\diam \bar g_k (Z) < \d _{k+1}/3$.
Let $g_{k+1}|_{\GB{k+1}} \equiv \bar g_k$. Observe, that $g_{k+1}$ is
continuous (see property G\ref{G2} on the page \pageref{G2}). Check
that $g_{k+1}$ satisfies conditions 1g--4g.
Indeed, 1g and 2g are met by construction. Further, we have
$\diam g_{k+1}(\GB{k+1}) < \d_ {k+1}/3$, therefore, since $y_{\al} \in
g_{k+1}(U_{\al})$ for every ${\al} \in \A$, condition 4g is also met.
Finally, our choice of $\d$ and $\{ \d_k \}$ coupled with property 3g for
$g_k$ and just checked property  4g for $g_{k+1}$ (as well as the choice
of $y_{\al}$) yields the property 3g for $g_{k+1}$.

Since $F_{n+1} = X$ (see property  F\ref{f3} on the page \pageref{f3})
we complete the proof letting $g \equiv g_{n+1}$.
\end{proof}

\begin{lem}\label{l2}
Let $B$ be a closed subset of $\AE $-compactum $Y$, such that $B \in
\AE ([L])$.
Then

{
\renewcommand{\theenumi}{{\bf \alph{enumi}}}
\begin{enumerate}
\item \label{a} For any $\m > 0$ there exists $\n = \n (\m ) > 0$
such that for any $X$ with $\ed (X) \le [L]$, any $A \subseteq X$ closed
in $X$ and any map $f: A \to O_{\n}(B)$ there exists $\bar f: X \to
O_{\m}(B)$ extending $f$.
\item \label{b} For any $\e > 0$ there exists $\d = \d (\e ) >0$
such that for any $\m >0$ there exists $\n = \n (\e , \m )$ with the
following property:

for every $X$ with $\ed (X)\le [L]$, any $A \subseteq X$ closed in $X$
and any $f: A\to O_{\n}(B)$ such that $\diam f(A) < \d$
there exists\\ $\bar f: A\to O_{\m}(B)$ such that $\diam f(X) < \e$.
\end{enumerate}
}
\end{lem}

\begin{proof}
{\bf a.} Since $Y$ is $\AE $-compactum, for $\m >0$ pick $\x = \x (\m )$ such that pair $(\m , \x )$
satisfies conditions of Corollary~\ref{c1} for space $Y$. Further, for $\x > 0$
choose $\n$ such that pair $(\x , \n )$ meets conditions of Lemma
\ref{l1}. Check that pair $(\m , \n)$ satisfies condition {\bf a}.

Consider $f: A \to O_{\n}(B)$. By our choice of $\n$ there exists $g: A
\to B$ such that $\dist (f,g) < \x$. Since $B \in \AE (X)$ there exists
extension $\bar g: X \to B$. Therefore, by our choice of $\x$, there
exists $\bar f: X \to Y$ such that $\dist (\bar f,
\bar g) < \m$.

The last fact implies that $f(X) \subseteq O_{\m}(B)$.

{\bf b.} Consider $\e ' = \e / 4$. For $\e '$ pick $\d ' = \d (\e ') > 0 $
such that pair $(\e ' , \d ')$ satisfies conditions of Proposition
\ref{hext} for space $B$. Let $\d = \d '/4$. Observe, that $\d = \d (\e )$.

Further, let $\l = \min (\m , \frac{\e}{3})$. For $\l$ find $\x > 0$
as in Corollary~\ref{c1}, applied to space $Y$. We may assume that
$\x < \d$. For $\x > 0$ pick $\n = \n (\x ) > 0$ as in Lemma~\ref{l1}.

Check, that $\d$ and $\n$ satisfy our requirments.

Consider $f: A \to O_{\n}(B)$ such that $\diam f(A) < \d$. By the choice
of $\n$ there exists $g: A \to B$, which is $\x$-close to $f$ and hence
$\d$-close to $f$. This fact implies, that $\diam g(A) < 3\d = \d '$,
which, in turn, implies by the choice of $\d '$, that $g$ has an
extension $\bar g: X \to B$ such that $\diam \bar g (X) < \e$.

Since $g$ and $f$ are $\x$-close, by our choice of $\x$ we may now
conclude that $f$ has an extension $\bar f: X \to Y$ such that $\bar f$
is $\l$-close to $\bar g$.
Finally, by the choice of $\l$, $\diam \bar f (X) < \e$ and
$\bar f (X) \subseteq O_{\m} (B)$.
\end{proof}

\section{Proof of Main Lemma}

\begin{proof}
Let we are given a  $\d >0$
and $\V \in \cov (X)$. We say, that cover $\U = \{ U_{\al} :
\al \in \A \}$ and sequences $\{\U _k \in \cov (X) : k=1 \ldots n+1\}$,
$\{ x(k,\al ) : \al \in \A, k = 1 \ldots n+1 \}$
form {\it canonical system with respect to} $F$ {\it and}
$\d$,
if the following conditions are satisfied (we use notation of
Definitions~\ref{can},~\ref{canF}):
{
\begin{enumerate}
\item \label{u1} $\U _{1} \equiv \V$.
\item \label{u2} $\U_{k+1}$ is canonical refinement of $\U_{k}$
with respect to $\d$ and $F$.
\item \label{u3} $\U$ is canonical refinement of $\U _{n+1}$.
\item \label{u4} $\St _{\U} U_{\al} \subseteq U_{x(n+1,\al )} \in \U
_{n+1}$ such that
$x(n+1, \al ) = x (U_{x(n+1,\al )})$.
\item \label{u5} $\St _{\U _{k+1}} U_{x(k+1,\al )} \subseteq U_{x(k,\al
)} \in \U _k$ such that
$x(k, \al ) = x (U_{x(k,\al )})$.
\end{enumerate}
}

Note, that canonical system exists for each $\V \in \cov (X)$.
Note also, that some of $\{ x(k, \al ) \}$ may coincide.

Finally, observe, that since $\{ F(x) : x \in X \}$ is $LC^{[L]}$
collection, we may assume without loss of generality
that $\d$ and $\n$ which Lemmas~\ref{l1},~\ref{l2} provide us with for
every $F(x)$ do not depend on $x$.

{\bf a.} Fix $\m >0$. Construct sequence $\{\d_k \}_{k=1}^{n+1}$
such that $\d _{n+1} = \m$ and for each $k$ pair $(\d _{k+1}/2,\d _k)$
satisfies conditions of Lemma~\ref{l2}.{\bf a} for any $F(x)$.
Let $\d = \d _1 /2$.  In addition, we may assume that $\d _k < \d _{k+
1}$ for every $k$.

Consider also a cover $\V = \{ U_x : x \in X\}$, where
$U_x = \{x' \in X: F(x) \subseteq  O_{\d} (F(x'))\}$.

Let $\U = \{ U_{\al} :\al \in \A \}$, $\{\U _k \in \cov (X) : k=1 \ldots n+1\}$,
$\{ x(k,\al ) : \al \in \A, k = 1 \ldots n+1 \}$ be canonical system for
$\d$  and $\V$.

For each $\al \in \A$ pick $y_{\al} \in F(x(1, \al))$.

Finally, consider sets $\{ F_k \}$ and $\{ \GB{k}  \}$, constructed with
respect to $\U$ (see Section~\ref{cov}).

Now we construct map $g$.

Using induction by $k = 1 \ldots n+1$ construct a sequence
of maps $\{ g_k \} ^{n+1}_{k=1}$ such that:
{\renewcommand{\theenumi}{\roman{enumi}}
\begin{enumerate}
\item \label{1f} $g_k: F_k \to Y$.
\item \label{2f} $g_{k+1}|_{F_k} = g_k$.
\item \label{3f} $g_k$ is $\d _k$-close to $F|_{F_k}$.
\item \label{4f} For each $\B \subseteq \A$ such that
$|\B | = i \le k$ there exists $\al = \al (\B ) \in \B$ having property
$g_k(\GB{i}) \subseteq O_{\d _i /2} ( F(x(i, \al (\B ))))$.
\end{enumerate}}
For $k = 1$ define $g_1$ letting $g_1|_{G^{ \{ \al \} } _1 } \equiv y(\al)$.
Observe, that $g_1$ is defined correctly and continuous on $F_1$
(see properties G\ref{G1}--G\ref{G4} on the page \pageref{G1}).
By properties~\ref{u1}--\ref{u4} of canonical system and by the choice
of $\{ y_{\al }\}$, $g_1$ satisfies requirements
\ref{1f}--\ref{4f}.

Assuming that $g_k$ has been already constructed, let us construct
$g_{k+1}$.

To accomplish this it is enough (see G\ref{G3} on the page  \pageref{G3})
to define $g_{k+1}$ on each non-empty $G^{\B }_{k+1}$
for each $\B \subseteq \A$ such that $|\B | = k+1$.

Fix $\B$ such that $|\B | = k+1$. Consider $Z = \overline{\GB{k+1}}$.
Let $Z' = Z \bigcap F_k$.
Obviously, $Z'$ is closed and non-empty subset of $Z$.
The idea is to define $g_{k+1}$ on $G^{\B }_{k+1}$
extending $g_k$ from $Z'$ over $Z$.

For each $\B ' \subsetneq \B$ consider $\al (\B ')$ which exists by
propery~\ref{4f} for map $g_k$.  Consider $U_{x(k+1,\al (\B ' ))}$.

By properties~\ref{u4} and~\ref{u5} of canonical system we have:
$$
U_{\al (\B ')} \subseteq U_{x(k+1, \al (\B '))}
$$
$(*)$ \phantom{**************************}  and
$$
\St _{\U _{k+1}} U _{x(k+1, \al (\B '))} \subseteq  U_{x(k, \al (\B '))}
\subseteq U_{x(|\B '|, \al ( \B '))}
$$

Since
$\bigcap \{ U_{x(k+1, \al (\B '))} : \B ' \subsetneq \B \} \supseteq
\GB{k+1} \ne \emp$ from $(*)$ we can conclude that there exists
$\B '' \subsetneq \B$ with the following property:
\begin{center}$(**)$\end{center}
$$
U_{x(k+1, \al (\B ''))} \subseteq \bigcap\{ U_{x(k+1, \al )}
: \al \in \B \} \subseteq \bigcap\{ U_{x(|\B '|, \al(\B ') )} :
\B ' \subsetneq \B \}
$$

Define $\al (\B ) = \al (\B '')$.

Property $(**)$ coupled with property~\ref{ex2} of Definition~\ref{canF}
implies that for any $\B ' \subsetneq \B$ we have
$F(x(k+1, \al (\B '))) \subseteq O_{\d} (F(x(k+1, \al (\B ))))$. Last
inclusion and property~\ref{4f} of $f_k$ (as well as our
choice of the sequence $\{\d _l\}$)  yields the following
chain of inclusions for each $\B ' \subsetneq \B$, $|\B '| = i \le k$ :

$g_k(\GBP{|\B '|}) \subseteq O_{\d _i/2 +\d} (F (x (k+1, \al (\B ))))
\subseteq O_{\d _k/2 +\d} (F (x (k+1, \al (\B ))))
\subseteq O_{\d _k} (F (x (k+1, \al (\B ))))$.

Hence (see property G\ref{G4} on the page \pageref{G4})
$g_k(Z') \subseteq O_{\d _k} (F (x (k+1, \al (\B ))))$.

The last fact and our choice of sequence $\{\d _l\}$ allow us to
extend $g_k$ to $g_{k+1}$ over $Z$ such that
$$
g_{k+1}(Z) \subseteq O_{\d _{k+1}/2} (F (x (k+1, \al (\B ))))
\leqno{(***)}
$$

Observe, that $g_{k+1}$ is correctly defined and continuous on $F_{k+1}$
by the properties  G\ref{G2} and G\ref{G4} on the page \pageref{G2}.
Let us check that $g_{k+1}$ satisfies conditions~\ref{1f}--\ref{4f}.

Indeed, conditions~\ref{1f} and~\ref{2f} are met by construction.

Further, since $\GB{k+1} \subseteq Z$, condition~\ref{4f} follows from
$(***)$.
Finally, since $\GB{k+1} \subseteq U_{x(k+1, \al (\B ))}$,  from
property~\ref{ex2} of Definition~\ref{canF}
applied to $\U_{k+1}$ we have $\dist (F|_{\GB{k+1}}, g_{k+1}|_{\GB{k+1}})
< \d _{k+1}/2 + \d < \d _{k+1} < \m$, which shows that property~\ref{3f}
is also met.

Since $F_{n+1} = X$ (see property  F\ref{f3} on the page \pageref{f3})
we complete the proof letting $g \equiv g_{n+1}$.


{\bf b.} Fix $\m , \e >0$. Construct sequences $\{\d _k \}_{k=1}^{n+1}$,
$\{\n _k \}_{k=1}^{n+1}$
such that $\d _{n+1} = \e$, $\n _{n+1} = \m$
and for each $k$ we have $\d _k = \d (\d _{k+1}/6)$,
$\n _k = \n (\d_{k+1}/6, \n _{k+1}/2)$
in the sense of Lemma~\ref{l2}.{\bf b} applied to  $F(x)$ (for any $x
\in X$).

Let $\d = \d _1 /12$ and $\n = \n _1 /2$.  In addition, we may assume
that $\d _k < \d _{k+1} /2$ and $\n _k < \n _{k+1}$ for every $k$.

Suppose that we are given a map $f: X \to Y$ such that $f$ is $\d$-close
to $F$.

Consider a cover $\V = \{ V_x : x \in X\}$, where
$V_x = \{x' \in X: F(x) \subseteq  O_{\n} (F(x'))\} \bigcap
f^{-1}(O_{\d}f(x))$.

Let $\U = \{ U_{\al} :\al \in \A \}$, $\{\U _k \in \cov (X) : k=1 \ldots n+1\}$,
$\{ x(k,\al ) : \al \in \A, k = 1 \ldots n+1 \}$ be canonical system for
$\n$  and $\V$.

Since $\dist (f,F) < \d$, for each $\al \in \A$ we can
pick $y_{\al} \in F(x(1, \al))$ such that $\dist (y_{\al}, f(x(1,\al
))) < \d$.

Finally, consider sets $\{ F_k \}$ and $\{ \GB{k}  \}$, constructed with
respect to $\U$ (see Section~\ref{cov}).

Now we construct map $g$.

As in part {\bf a}, using induction by $k = 1 \ldots n+1$
construct a sequence
of maps $\{ g_k \} ^{n+1}_{k=1}$ such that:
{\renewcommand{\theenumi}{\roman{enumi}}
\begin{enumerate}
\item \label{1b} $g_k: F_k \to Y$.
\item \label{2b} $g_{k+1}|_{F_k} = f_k$.
\item \label{3b} $g_k$ is $\n _k$-close to $F|_{F_k}$.
\item \label{4b} For each $\B \subseteq \A$ such that
$|\B | = i \le k$ there exists $\al = \al (\B ) \in \B$ having property
$g_k(\GB{i}) \subseteq O_{\n _i /2} ( F(x(i, \al (\B ))))$.
\item \label{5b} For each $\al$,
$g_k (F_k \bigcap U_{\al}) \subseteq O_{\d _k/3}(y_{\al})$
\end{enumerate}}
For $k = 1$ define $g_1$ letting $g_1|_{G^{ \{ \al \}} _1}  \equiv y(\al)$.
Observe, that $g_1$ is defined correctly and continuous on $F_1$
(see properties~G\ref{G1}--G\ref{G4} on the page \pageref{G1}).
By properties~\ref{u1}--\ref{u4} of canonical system and by the choice
of $\{ y_{\al }\}$, $g_1$ satisfies requirements
\ref{1b}--\ref{5b}.

Assuming that $g_k$ has been already constructed, let us construct
$g_{k+1}$.

As before in the proof of {\bf a}, to accomplish this it is enough (see G\ref{G3} on the page  \pageref{G3})
to define $g_{k+1}$ on each non-empty $G^{\B }_{k+1}$
for each $\B \subseteq \A$ such that $|\B | = k+1$.

Fix $\B$ such that $|\B | = k+1$. Consider $Z = \overline{\GB{k+1}}$.
Let $Z' = Z \bigcap F_k$.
Obviously, $Z'$ is closed and non-empty subset of $Z$.
Again, the idea is to define $g_{k+1}$ on $G^{\B }_{k+1}$
extending $g_k$ from $Z'$ over $Z$.

Using the same arguments as in proof of {\bf a}, one can show, that
$$
g_k(Z') \subseteq O_{\n _k} (F (x (k+1, \al (\B ))))\eqno{(')}
$$

Let us show, that
$$
\diam g_k (Z') < \d _k\eqno{('')}
$$
Indeed, since $\bigcap\limits_{\al \in \B} U_{\al} \supseteq \GB{k+1}
\ne \emp$, we have $\bigcap\limits_{\al \in \B} U_{x(1,\al )} \ne \emp$.
Therefore $\dist (f(x(1, \al )), f(x(1, \be ))) < 2\d$ for any $\al
, \be \in \B$.
Further, by construction we have $\dist (y_{\al}, f(x(1,\al ))) < \d$
for any $\al \in \B$. These inequalities coupled with property~\ref{5b}
and the fact that $Z' \subseteq \bigcup\limits_{\al \in \B} U_{\al}$
yield $\diam g_k (Z') <2\d + 2\d + 2(\d_k/3) = 4\d +2(\d _k/3) < \d _k$.
Property $('')$ is checked.

 From properties $(')$ and $('')$ we can conclude according to our
choice of sequences $\{\d _i \}$ and $\{\n _i\}$ that
$g _k$ can be extended over $Z$ to $g_{k+1}: Z \to Y$ such that
$$
\diam g_{k+1} (Z) < \d _{k+1}/6\eqno{(''')}
$$

Using the same arguments as in proof of {\bf a} one can show
that $g_{k+1}$ is continuous map satisfying properties~\ref{1b}--\ref{4b}.

Let us check that property~\ref{5b} is also met.

Since $Z' \bigcap U_{\al} \ne \emp$, for each $\al \in \B$,
from property $(''')$ and property~\ref{5b} applied to
map $g_k$ we have:
$$
g_{k+1}(F_{k+1} \bigcap U_{\al}) \subseteq O_{\d _{k+1}/6 + \d _k/3}
(y_{\al}) \subseteq O_{\d _{k+1}/3}(y_{\al})
$$
and condition~\ref{5b} is checked.

Finally, let $g \equiv g_{n+1}$.

Obviously, $g$ is $\m$-close to $F$.
Check, that $g$ is $\e$-close to $f$.

Property~\ref{5b} of $g_{n+1} = g$ implies that
$g_{n+1} (U_{\al}) \subseteq O_{\d _{n+1}/3} (y _{\al}) = O_{\e /3} (y
_{\al})$.
Since $U_{\al} \subseteq U_{x(1, \al )} \subseteq f^{-1}(O_\d f(x(1,
\al )))$, we have $f(U_{\al})  \subseteq O_{\d}f(x(1,\al ))$.

These inclusions coupled with inequality $\dist(f(x(1,\al )), y_{\al}) <
\d$ imply that
$\dist (g|_{U_{\al}}, f|_{U_{\al}}) < 2\d + \e/3 <\e$ for each $\al $
and consequently
$g$ is $\e$-close to $f$.
\end{proof}

The author is grateful to A.~C.~Chigogidze for attension to this work
and useful discussions.



\begin{thebibliography}{99}

\bibitem{mchl}
E.~Michael, {\em Continuous selections II}, Ann. Math. {\bf 63} (1956),
562-580.


\bibitem{ch}
A.~Chigogidze, {\em Infinite dimensional topology and shape theory}, to
appear in: "Handbook of Geometric Topology" edited by R.~Daverman and
R.~B.~Sher), North Holland, Amsterdam, 1999.


\bibitem{eng}
R.~Engelking, {\em General Topology}, PWN, W rszawa, 1977.


\end{thebibliography}
\end{document}